\documentclass[12pt]{amsart}
\usepackage{amsmath,amssymb}
\usepackage{amsfonts}
\usepackage{amsthm}
\usepackage{latexsym}
\usepackage{graphicx}
\def\p{\partial}

\def\R{\mathbb{R}}

\def\vv<#1>{\langle#1\rangle}

\def\1{\mathbf{1}}

\def\XXint#1#2{\setbox0=\hbox{$#1{#2}{\int}$}{#2}\kern-.5\wd0 }

\def\XXint#1#2#3{{\setbox0=\hbox{$#1{#2#3}{\int}$}
     \vcenter{\hbox{$#2#3$}}\kern-.5\wd0}}

\def\vv<#1>{{\left\langle#1\right\rangle}}

\def\CD{{\rm CD}}
\def\Deg{{\rm Deg}}

\newtheorem{thm}{Theorem}[section]

\newtheorem{prop}{Proposition}[section]
\newtheorem{cor}{Corollary}[section]
\theoremstyle{definition}
\newtheorem{defn}{Definition}[section]
\theoremstyle{remark}

\newtheorem{rem}{Remark}[section]

\numberwithin{equation}{section}

\begin{document}
\title{Comparisons of Dirichlet, Neumann and Laplacian eigenvalues on graphs and their applications}

\author{Yongjie Shi$^1$}
\address{Department of Mathematics, Shantou University, Shantou, Guangdong, 515063, China}
\email{yjshi@stu.edu.cn}
\author{Chengjie Yu$^2$}
\address{Department of Mathematics, Shantou University, Shantou, Guangdong, 515063, China}
\email{cjyu@stu.edu.cn}
\thanks{$^1$Research partially supported by NNSF of China with contract no. 11701355. }
\thanks{$^2$ Research partially supported by GDNSF with contract no. 2021A1515010264 and NNSF of China with contract no. 11571215.}
\renewcommand{\subjclassname}{%
  \textup{2010} Mathematics Subject Classification}
\date{}
\begin{abstract}
 In this paper, we obtain some comparisons of the Dirichlet, Neumann and Laplacian eigenvalues on graphs. We also discuss their rigidities and some of their applications including some Lichnerowicz-type, Fiedler-type and Friedman-type estimates for Dirichlet  eigenvalues and Neumann  eigenvalues. The comparisons on Neumann eigenvalues can be translated to comparisons on Steklov eigenvalues in our setting. So, some of the results can be viewed as extensions for parts of the works of \cite{HHW} by Hua-Huang-Wang, and parts of our previous works \cite{SY,SY2}.
\end{abstract}
\maketitle\markboth{Shi \& Yu}{Comparison of  eigenvalues on graphs}
\section{Introduction}
Let $(M^n,g)$ be a closed Riemannian manifold with Ricci curvature bounded from below by a positive constant $K$. Then, the well-known Lichnerowicz estimate \cite{LI} tells us that the first positive Laplacian eigenvalue of $(M^n, g)$ is no less than $\frac{n K}{n-1}$. This estimate was later extended to compact Riemannian manifolds with boundary by Reilly \cite{RE}. In recent years, Lichnerowicz estimate was extended to graphs in \cite{BC, KKR, LLY}. So, it is a natural problem to extend Reilly's Lichnerowicz estimate to graphs with boundary.

On the other hand, in  the recent works \cite{SY,SY2}, the authors obtained Lichnerowicz estimates for Steklov eigenvalues on graphs which may be viewed as an extension of the works of Escobar \cite{ES} and Xia-Xiong \cite{XX} into  discrete setting, by using a comparison of Steklov eigenvalues and Laplacian eigenvalues on graphs that was also mentioned in \cite{HHW2} for graphs with normalized weights.  It seems that such kinds of eigenvalue comparisons make a major difference of spectral theory on graphs with that on Riemannian manifolds. In this paper, motivated by our previous works, by further exploring comparisons of Dirichlet, Neumann  and Laplacian eigenvalues on graphs, we obtain Lichenerowicz-type estimates for Dirichlet and Neumann eigenvalues on graphs  extending the classical results of Reilly \cite{RE} into  discrete setting.

Firstly, we have the following comparison of Neumann eigenvalues and Laplacian eigenvalues. We use $\nu_i$ and $\mu_i$ to denote the $i^{\rm th}$ Neumann eigenvalue and $i^{\rm th}$ Laplacian eigenvalue of a weighted graph. For detailed definitions and other notations, see Section 2.
\begin{thm}\label{thm-Neumann-Laplacian}
Let $(G,B,m,w)$ be a  weighted finite graph with boundary. Then,
\begin{equation}\label{eq-Neumann-Laplacian}
\nu_i\geq\mu_i
\end{equation}
for any $i=1,2,\cdots, |\Omega|$. If $\mu_i=\nu_i>0$ for some $i=2,\cdots,|\Omega|$, then there is a function $u\in \R^V$ such that $u|_B=\frac{\p u}{\p n}=0$ and $u$ is simultaneously a Laplacian eigenfunction and a Neumann Laplacian eigenfunction of the eigenvalue $\nu_i=\mu_i$.
\end{thm}
We next have the following rigidity for \eqref{eq-Neumann-Laplacian}. For a weighted graph $(G,m,w)$, and $S\subset V(G)$, we use $\mu_i(S)$ to denote the $i^{\rm th}$ Laplacian eigenvalue of the induced subgraph on $S$ inheriting vertex-measure and edge-weight from $G$. For detailed definitions and other notations, see Section 2.
\begin{thm}\label{thm-Neu-Lap-rigidity}
Let $(G,B,m,w)$ be a connected weighted finite graph with boundary. Then, the equalities of
\eqref{eq-Neumann-Laplacian}
hold for any $i=1,2,\cdots, |\Omega|$ if and only if
the following two statements are both true.

(1) There is a nonnegative function $\rho\in \R^B$, such that $w_{xy}=\rho_xm_xm_y$ for any $x\in B$ and $y\in \Omega$. In particular, each boundary vertex is either adjacent to every interior vertex or adjacent to no interior vertex.

(2) When $\rho$ is constant,
\begin{equation}\label{eq-com-Neu-Lap-rho-const}
\mu_{|\Omega|}(\Omega)\leq \min\left\{\rho V_\Omega, \mu_2(B)+\rho(V_\Omega-V_B)\right\}.
\end{equation}
When $\rho$ is nonconstant,
\begin{equation*}
\mu_{|\Omega|}(\Omega)<\frac{V_\Omega}{V_B}\Deg_b,
\end{equation*}
and for any $f\in \R^B$ with $\vv<f,1>_B=0$,
\begin{equation*}
\vv<df,df>_B+V_\Omega\vv<\rho f,f>_B-\left(\mu_{|\Omega|}(\Omega)+\Deg_b\right)\vv<f,f>_B\geq \frac{V_\Omega V_G}{V_\Omega\Deg_b-V_B\mu_{|\Omega|}(\Omega)}\vv<\rho,f>_B^2.
\end{equation*}
Here
$$\Deg_b:=\vv<\rho,1>_B=\Deg_B(y)$$ for any $y\in \Omega$, $V_B=\sum_{x\in B}m_x$, $V_\Omega=\sum_{y\in\Omega}m_y$ and $V_G=V_\Omega+V_B.$
\end{thm}

By Theorem \ref{thm-Neu-Lap-rigidity}, one can construct many nontrivial examples such that equality of \eqref{eq-Neumann-Laplacian} holds for $i=1,2,\cdots, |\Omega|$. For example, fix a graph $(G,B)$ with boundary  containing the complete bipartite graph $K_{B,\Omega}$. Set the weights of boundary edges and measures of interior vertices such that $V_\Omega>V_B$, and $w_{xy}=\rho m_xm_y$ for any $x\in B$ and $y\in \Omega$. Finally, set the weights  of interior edges small enough to make
 $$\mu_{|\Omega|}(\Omega)\leq \rho(V_\Omega-V_B).$$
Then, by Theorem \ref{thm-Neu-Lap-rigidity}, equalities of \eqref{eq-Neumann-Laplacian} hold for all $i=1,2,\cdots, |\Omega|$ on  the graph constructed.

As an application of Theorem \ref{thm-Neu-Lap-rigidity}, we have the following rigidities of \eqref{eq-Neumann-Laplacian} for graphs equipped with the unit weight or a normalized weight.
\begin{cor}\label{cor-rigidity-unit-weight}
Let $(G,B)$ be a connected finite graph with boundary equipped with the unit weight such that each boundary vertex is adjacent to some interior vertex. Then, the equalities of
\eqref{eq-Neumann-Laplacian}
hold for any $i=1,2,\cdots, |\Omega|$ if and only if
\begin{enumerate}
\item every boundary vertex is adjacent to any interior vertex, and
\item $\mu_{|\Omega|}(\Omega)\leq \mu_2(B)+|\Omega|-|B|$ when $|B|\geq 2$, and $\mu_{|\Omega|}(\Omega)\leq |\Omega|$ when $|B|=1$.
 \end{enumerate}
In particular, if every boundary vertex is adjacent to any interior vertex and  $$2\Delta(G|_\Omega)\leq|\Omega|-|B|,$$ then the equalities of \eqref{eq-Neumann-Laplacian} hold for $i=1,2,\cdots,|\Omega|$. Here $\Delta(G|_\Omega)$ is the maximum degree of $G|_\Omega$.
\end{cor}
\begin{cor}\label{cor-rigidity-normalized-weight}
Let $(G,B,m,w)$ be a connected finite graph with boundary equipped with a normalized weight such that $E(B,B)=\emptyset$. Then, the equalities of
\eqref{eq-Neumann-Laplacian}
hold for any $i=1,2,\cdots, |\Omega|$ if and only if
\begin{enumerate}
\item for any $x\in B$ and $y\in \Omega$, $w_{xy}=\frac{1}{V_\Omega}m_xm_y$, and
\item $\mu_{|\Omega|}(\Omega)\leq \frac{V_\Omega-V_B}{V_\Omega}$ when $|B|\geq 2$, and $\mu_{|\Omega|}(\Omega)\leq 1$ when $|B|=1$.
 \end{enumerate}
\end{cor}

Secondly, we have the following comparisons of Dirichlet eigenvalues, Neumann eigenvalues and Laplacian eigenvalues of the interior subgraph. We use $\lambda_i$ to denote the $i^{\rm th}$ Dirichlet eigenvalue of a weighted graph. For details of definitions and other notations, see Section 2.
\begin{thm}\label{thm-D-L}
Let $(G,B,m,w)$ be a weighted finite graph with boundary. Then,
\begin{enumerate}
\item for any $i=1,2,\cdots,|\Omega|$ and $j=0,1,\cdots,|\Omega|-i$,
\begin{equation*}
\lambda_i\leq \mu_{i+j}(\Omega)+\Deg_B^{(|\Omega|-i)}
\end{equation*}
with equality for some pair $i, j$ if and only if there is a nonzero function $u\in \R^\Omega$ such that $\Delta^D u=\lambda_i u$, $\Delta_\Omega u=\mu_{i+j}(\Omega)u$ and $\Deg_B \cdot u=\Deg_B^{(|\Omega|-i)}\cdot u$;
\item for any $i=1,2,\cdots,|\Omega|$ and $j=1,2,\cdots, i$,
\begin{equation*}
\lambda_i\geq \mu_{i-j+1}(\Omega)+\Deg_B^{(j)}
\end{equation*}
with equality for some pair $i, j$ if and only if there is a nonzero function $u\in \R^\Omega$ such that $\Delta^D u=\lambda_i u$, $\Delta_\Omega u=\mu_{i-j+1}(\Omega)u$ and $\Deg_B\cdot u=\Deg_B^{(j)}\cdot u$.
\end{enumerate}
Here
$$\Deg_B^{(1)}\leq \Deg_B^{(2)}\leq \cdots\leq \Deg_B^{|\Omega|}$$
is the rearrangement of $\Deg_B(x)$ with $x$ going through $\Omega$.
\end{thm}
\begin{thm}\label{thm-N-L}
Let $(G,B,m,w)$ be a weighted finite graph with boundary. Then
\begin{equation*}
\nu_i\geq \mu_i(\Omega)
\end{equation*}
for $i=1,2,\cdots,|\Omega|$, with equalities for $i=1,2,\cdots,|\Omega|$ if and only if each connected component of $N_G(B)$ contains at most one interior vertex. Here
$$N_G(B)=B\cup\{y\in V(G) \ |\ y\sim x\mbox{ for some }x\in B\}.$$
\end{thm}
\begin{thm}\label{thm-N-D}
Let $(G,B,m,w)$ be a  weighted finite graph with boundary. Then,
$$\lambda_i\geq \nu_i$$
for $i=1,2,\cdots,|\Omega|$.
\end{thm}
Combining Theorem \ref{thm-Neumann-Laplacian} and Theorem \ref{thm-N-D}, we have the following comparison of Dirichlet eigenvalues and Laplacian eigenvalues.
\begin{cor}\label{cor-D-L}
Let $(G,B,m,w)$ be a  weighted finite graph with boundary. Then,
$$\lambda_i\geq \mu_i$$
for $i=1,2,\cdots,|\Omega|$.
\end{cor}
Combining (1) of Theorem \ref{thm-D-L} and Corollary \ref{cor-D-L}, we have the following comparison of the Laplacian eigenvalues and the Laplacian eigenvalues of the interior subgraph.
\begin{cor}\label{cor-L-IL}
Let $(G,B,m,w)$ be a weighted finite graph with boundary. Then,
$$\mu_i\leq\mu_{i+j}(\Omega)+\Deg_B^{(|\Omega|-j)}$$
for $i=1,2,\cdots,|\Omega|$ and $0,1,\cdots,|\Omega|-i$.
\end{cor}
Note that $\mu_i$ is independent of the choices of $B$. Then, we have the following upper bound estimate for $\mu_i$ by Corollary \ref{cor-L-IL} and noting that $\mu_i(\Omega)=0$ if $G|_\Omega$ has at least $i$ connected components. It can be viewed as an extension of the vertex-connectivity upper bounds for Laplacian eigenvalues by Fiedler \cite{FI}.
\begin{cor}\label{cor-upper}
Let $(G,m,w)$ be a weighted finite graph. Then, for $i=1,2,\cdots,|V|$,
$$\mu_i\leq \min\left\{\Deg_S^{(|S^c|-j)}\ \bigg|\  S\subset V,\ 0\leq j\leq |S^c|-i,\mbox{ and }\pi_0(G-S)\geq i+j \right\}.$$
In particular, for $i=1,2,\cdots,|V|$,
$$\mu_i\leq \min\left\{\max_{x\in S^c}\Deg_S(x)\ \bigg|\  S\subset V\mbox{ such that }\pi_0(G-S)\geq i \right\}.$$
Especially, if $G$ is equipped with the unit weight, then
$$\mu_i\leq\min \left\{|S|\ |\ S\subset V\mbox{ such that }\pi_0(G-S)\geq i \right\},$$
for $i=1,2,\cdots,|V|$. Here $\pi_0(G-S)$ is the number of connected components of $G-S$, and
$$\Deg_S^{(1)}\leq \Deg_S^{(2)}\leq \cdots\leq \Deg_S^{|S^c|}$$
is the rearrangement of $\Deg_S(x)$ with $x$ going through $S^c$.
\end{cor}

Note that the equalities in Theorem \ref{thm-N-D} can not hold. For example, when $G$ is connected, $\nu_1=0$, but $\lambda_1>0$. This means that the comparison in Theorem \ref{thm-N-D} is not sharp. By comparing the Neumann Laplacian operator and the Dirichlet Laplacian operator, we have the following sharp comparison of Neumann eigenvalues and Dirichlet eigenvalues.
\begin{thm}\label{thm-N-D-sharp}
Let $(G,B,m,w)$ be a  weighted finite graph with boundary such that each connected component of $G$ contains some interior vertices. Then
\begin{equation}\label{eq-N-D-sharp}
\nu_i+\lambda_{\min}(A_BL_B^{-1}A_\Omega)\leq\lambda_i\leq \nu_i+ \lambda_{\max}(A_BL_B^{-1}A_\Omega)
\end{equation}
for $i=1,2,\cdots, |\Omega|$, with equalities for $i=1,2,\cdots,|\Omega|$ if and only if $A_BL_B^{-1}A_\Omega:\R^\Omega\to \R^\Omega$ is a scalar operator.
\end{thm}
Here $A_B$ and $A_\Omega$ are the weighted adjacent operators relative to $B$ and $\Omega$ respectively (see Section 2 for details), and $L_B:\R^B\to \R^B$ is the Dirichlet Laplacian operator with Dirichlet boundary value on $\Omega$ (see Section 3 for details). Moreover $\lambda_{\min}(T)$ and $\lambda_{\max}(T)$ denote the minimal eigenvalue and maximal eigenvalue of the linear transformation $T$ respectively.

Because it is hard to get explicit expression of $L_B^{-1}$ for general weighted graphs with boundary, we are not able to obtain a geometric characterization of graphs that the equalities of \eqref{eq-N-D-sharp} all hold. However, when the boundary of the graph is simple in the sense that $E(B,B)=\emptyset$, we have an explicit expression for $L_B^{-1}$. With the help of this expression, we are able to get a more geometric characterization of graphs that the equalities of \eqref{eq-N-D-sharp} all hold.
\begin{thm}\label{thm-N-D-sharp-1}
Let $(G,B,m,w)$ be a connected weighted finite graph with boundary such that $E(B,B)=\emptyset$. Then
\begin{equation}\label{eq-com-Neu-Diri}
\nu_i+\lambda_{\min}(A_B\Deg^{-1}A_\Omega)\leq\lambda_i\leq\nu_i+\lambda_{\max}(A_B\Deg^{-1}A_\Omega)
\end{equation}
for any $i=1,2,\cdots, |\Omega|$. Moreover, the equalities of \eqref{eq-com-Neu-Diri} hold for $i=1,2,\cdots, |\Omega|$, if and only if
\begin{enumerate}
\item every boundary vertex is adjacent to only one interior vertex, and
\item the weighted degree $\Deg_B(x)$ relative to $B$
 is independent of $x\in \Omega$.
\end{enumerate}

\end{thm}
There are also comparisons between Dirichlet eigenvalues and Neumann eigenvalues, and comparisons between Neumann eigenvalues and Laplacian eigenvalues for the interior subgraph similar to Theorem \ref{thm-D-L}. See Theorem \ref{thm-N-D-sharp-2} and Theorem \ref{thm-N-IL}  in Section 5 for details.

Finally, by combining the comparison results above and some known lower bounds for Laplacian eigenvalues such as Lichnerowicz-type estimates in \cite{BC,LLY,KKR}, Fielder's estimate \cite{FI} and Friedman's estimate \cite{FR}, we can get some interesting estimates on Neumann eigenvalues and Dirichlet eigenvalues. For example, we have the following Lichnerowicz-type estimates for the first positive Neumann eigenvalue. For other estimates, see Section 6 for details.
 \begin{thm}\label{thm-Lich-BE-1}
Let $(G,B,m,w)$ be a connected weighted finite graph with boundary.
\begin{enumerate}
\item If $(G,m,w)$ satisfies the Bakry-\'Emery curvature-dimension condition $\CD(K,n)$ with $K>0$ and $n>1$, then $\nu_2\geq \frac{nK}{n-1}$.
\item If the Ollivier curvature of $(G,m,w)$ has a positive lower bound $\kappa$, then $\nu_2\geq \kappa$.
\end{enumerate}
\end{thm}

Note that in the Lichnerowicz estimates above, no analogue of boundary curvatures such as mean curvature or second fundamental form was involved. This is different with the Riemannian case.

At the end of this section, we would like to mention that we take the most general definition of graphs with boundary in this paper. That is, any pair $(G,B)$ with $B\subset V(G)$ is a graph with boundary (see \cite{YY1} for example).  For all the comparison results of eigenvalues, we don't need the assumption that $E(B,B)=\emptyset$ (see \cite{Pe,HHW} for examples). However, when considering rigidities of the comparisons of eigenvalues, we sometimes need this assumption to obtain more geometric characterizations. See Corollary \ref{cor-rigidity-unit-weight}, Corollary \ref{cor-rigidity-normalized-weight} and Theorem \ref{thm-N-D-sharp-1} for examples.

Moreover, because in this most general notion of graphs with boundary, one can interchange the roles of the interior and boundary. So, Neumann eigenvalues and Steklov eigenvalues can be viewed as dual of each other in the sense that
$$\sigma_i(G,S)=\nu_i(G,S^c),\ \forall S\subset V(G).$$
Therefore, the comparisons on Neumann eigenvalues in this paper also provide comparisons on Steklov eigenvalues on the most general graphs with boundary, and Theorem \ref{thm-Neumann-Laplacian}, Theorem \ref{thm-Neu-Lap-rigidity} and Theorem \ref{thm-Lich-BE-1} are actually extensions of our previous works \cite{SY,SY2} to the most general graphs with boundary. The comparisons on Neumann eigenvalues can then be translated as comparisons on Steklov eigenvalues. See Corollary \ref{cor-Steklov-0}, Corollary \ref{cor-Steklov-1}, Corollary \ref{cor-Steklov-2} and Corollary \ref{cor-Steklov-3} in Section 5 for examples.

The rest of the paper is organized as follows. In section 2, we introduce some notations and preliminaries. In Section 3, we obtain some simple but useful relations on the Dirichlet Laplacian operator, Neumann Laplacian operator and Laplacian operator of the interior subgraph. In Section 4, we obtain the comparison of Neumann eigenvalues and Laplacian eigenvalues and its rigidity. In Section 5, we obtain some comparisons on Dirichlet eigenvalues, Neumann eigenvalues and Laplacian eigenvalues of the interior subgraph. In Section 6, we obtain some applications of the comparisons of eigenvalues by combining them with some known lower bounds for Laplacian eigenvalues.

\section{Preliminaries}
In this section, we introduce some preliminaries and notations on analysis of graphs. For more details about analysis on graphs, see \cite{Gr} and \cite{Ch}.

Firstly, let's recall the notion of weighted graphs.
\begin{defn}
A weighted graph is a tripe $(G,m,w)$ where
\begin{enumerate}
\item $G$ is a simple graph;
\item $m:V(G)\to \R^+$ is called the vertex-measure;
\item $w:E(G)\to \R^+$ is called the edge-weight.
\end{enumerate}
Here $V(G)$ and $E(G)$ are the sets of vertices and edges of the graph $G$ respectively. We will simply write them as $V$ and $E$ respectively if no confusion was made.
\end{defn}
When $m\equiv 1$ and $w\equiv 1$, we say that $G$ is equipped with the unit weight.

For simplicity, we also view $w$ as a symmetric function on $V\times V$ by zero extension. Throughout this paper, graphs are assumed to be simple and finite.

For each $x\in V$, define the weighted degree $\Deg(x)$ at $x$ as
\begin{equation*}
\Deg(x)=\frac{1}{m_x}\sum_{y\in V}w_{xy}.
\end{equation*}
If $\Deg(x)=1$ for any $x\in V$, we say that $G$ is equipped with a normalized weight.

For any $S\subset V$, define the weighted degree $\Deg_S(x)$ relative to $S$ at $x$ as
$$\Deg_S(x)=\frac{1}{m_x}\sum_{y\in S}w_{xy}.$$
Moreover, define the weighted adjacency operator $A_S:\R^S\to \R^{S^c}$ relative to $S$ as
\begin{equation*}
A_Su(x)=\frac{1}{m_x}\sum_{y\in S}u(y)w_{xy},\ \forall u\in \R^S\mbox{ and }x\in S^c.
\end{equation*}
It is clear that
$$\Deg_S(x)=A_S(\mathbf 1)(x) \ \forall x\in S^c.$$

Equip $\R^V$ with the inner product:
\begin{equation*}
\vv<u,v>=\sum_{x\in V}u(x)v(x)m_x,\ \forall u, v\in \R^V.
\end{equation*}
For $S\subset V$, define
$$\vv<u,v>_S=\sum_{x\in S}u(x)v(x)m_x.$$
It is then clearly that  $A_S$ and $A_{S^c}$ are the adjoint operator of each other:
$$\vv<A_S u,v>_{S^c}=\vv<u,A_{S^c}v>_S,\ \forall u\in \R^S\mbox{ and }v\in \R^{S^c}.$$

Let  $A^1(G)$ be the space of skew-symmetric functions $\alpha$ on $V\times V$ such $$\alpha(x,y)=0\mbox{ when } x\not\sim y.$$
In fact, $\alpha$ is called a flow in \cite{Ba}. Equip $A^1(G)$ with the inner product:
\begin{equation*}
\vv<\alpha,\beta>=\sum_{\{x,y\}\in E}\alpha(x,y)\beta(x,y)w_{xy}=\frac12\sum_{x,y\in V}\alpha(x,y)\beta(x,y)w_{xy},\ \forall \alpha,\beta\in A^1(G).
\end{equation*}
For $S\subset E$, define
\begin{equation*}
\vv<\alpha,\beta>_{S}=\sum_{\{x,y\}\in S}\alpha(x,y)\beta(x,y)w_{xy}, \ \forall \alpha,\beta\in A^1(G).
\end{equation*}
For $S\subset V$, define
$$\vv<\alpha,\beta>_S=\vv<\alpha,\beta>_{E(S,S)}, \ \forall \alpha,\beta\in A^1(G).$$

For any $u\in \R^V$, define the differential $du\in A^1(G)$ of $u$ as
\begin{equation*}
du(x,y)=\left\{\begin{array}{ll}u(y)-u(x)&\{x,y\}\in E\\0&\mbox{otherwise.}\end{array}\right.
\end{equation*}
From this definition, we know that $du\equiv 0$ if and only if $u$ is constant on each connected component of $G$.

Let $d^*:A^1(G)\to \R^V$ be the adjoint operator of $d:\R^V\to A^1(G)$. The Laplacian operator on $\R^V$ is defined as
\begin{equation*}
\Delta=d^*d.
\end{equation*}
By direct computation,
\begin{equation*}
\Delta u(x)=\frac{1}{m_x}\sum_{y\in V}(u(x)-u(y))w_{xy}
\end{equation*}
for any $x\in V$. Moreover, by the definition of $\Delta$, it is clear that
\begin{equation}\label{eq-integration-by-part}
\vv<\Delta u,v>=\vv<du,dv>
\end{equation}
for any $u,v\in \R^V$. So $\Delta$ is a nonnegative self-adjoint operator on $\R^V$. Let
\begin{equation*}
0=\mu_1\leq\mu_2\leq \cdots \leq \mu_{|V|}
\end{equation*}
be the eigenvalues of $\Delta$ on $(G,m,w)$. It is not hard to see from \eqref{eq-integration-by-part} that the multiplicity of eigenvalue $0$ equals to the number of connected components of $G$. So, when $G$ is connected, $\mu_2>0$.

We will also denote $\Delta$ as $\Delta_G$, and denote $\mu_i$ as $\mu_i(G)$ or $\mu_i(G,m,w)$ when it is necessary. For any nonempty subset $S$ of  $V$, the Laplacian operator of the induced subgraph $G|_S$ (inherit the weights from $G$) is denoted as $\Delta_S$. The eigenvalues of $\Delta_S$ is denoted as
  $$0=\mu_1(S)\leq \mu_2(S)\leq\cdots\leq \mu_{|S|}(S).$$
For $i>|S|$, we take the convention that $\mu_i(S)=+\infty$.

Next, let's recall the notion of graphs with boundary.
\begin{defn}\label{def-graph-boundary}
 A pair $(G,B)$ is called  a graph with boundary if $G$ is a graph and $B$ is a nonempty subset of $V(G)$ which is called the boundary of the graph. The set $\Omega:=V\setminus B$ is called the interior of the graph. Edges in $E(B,\Omega)$ are called boundary edges.
\end{defn}

Let $(G,B,m,w)$ be a weighted graph with boundary. For any $u\in \R^V$ and $x\in B$, define the normal derivative of $u$ at $x$ as:
\begin{equation*}
\frac{\p u}{\p n}(x):=\frac{1}{m_x}\sum_{y\in V}(u(x)-u(y))w_{xy}=\Delta u(x).
\end{equation*}
Then, by \eqref{eq-integration-by-part}, one has  the following Green's formula:
\begin{equation*}
\vv<\Delta u,v>_\Omega=\vv<du,dv>-\vv<\frac{\p u}{\p n},v>_B.
\end{equation*}

We now introduce  Dirichlet eigenvalues, Neumann  eigenvalues  and Steklov eigenvalues on weighted graphs with boundary.

A real number $\lambda$ is called a Dirichlet  eigenvalue of $(G,B,m,w)$ if the following Dirichlet boundary value problem:
\begin{equation*}
\left\{\begin{array}{ll}\Delta u(x)=\lambda u(x)&x\in\Omega\\
u(x)=0&x\in B
\end{array}\right.
\end{equation*}
has a nontrivial solution. The corresponding operator for Dirichlet eigenvalues is the Dirichlet Laplacian operator introduced as follows.

Let $E_0:\R^\Omega\to\R^V$ be defined as
\begin{equation*}
E_0(u)(x)=\left\{\begin{array}{ll}u(x)& x\in\Omega\\
0&x\in B.
\end{array}\right.
\end{equation*}
Then, the Dirichlet  eigenvalues of $(G,B,m,w)$ are the eigenvalues of the operator $\Delta^D:\R^\Omega\to\R^\Omega$ with
\begin{equation*}
\Delta^D u=\Delta E_0(u)|_\Omega
\end{equation*}
which is called the Dirichlet Laplacian operator. By \eqref{eq-integration-by-part}, it is clear that
\begin{equation}\label{eq-int-by-part-Dirichlet}
\vv<\Delta^D u,v>_\Omega=\vv<\Delta E_0(u),E_0(v)>=\vv<dE_0(u),dE_0(v)>
\end{equation}
for any $u,v\in\R^\Omega$. So, $\Delta^D$ is a nonnegative self-adjoint operator on $\R^\Omega$. We denote its eigenvalues as
\begin{equation*}
0\leq \lambda_1\leq\lambda_2\leq \cdots\leq \lambda_{|\Omega|}.
\end{equation*}
We also denote $\Delta^D$ as $\Delta^D_{G,B}$, and denote $\lambda_i$ as
$\lambda_i(G,B)$ or $\lambda_i(G,B,m,w)$ when it is necessary. When $i>|\Omega|$, we take the convention that $\lambda_i=+\infty$.

By \eqref{eq-int-by-part-Dirichlet}, we have the following conclusion for the multiplicity of the eigenvalue $0$ for $\Delta^D$.
\begin{prop}\label{prop-ker-D}
Let $(G,B,m,w)$ be a weighted finite graph with boundary. Let $i(G,B)$ be the number of connected components of $G$ that contains no boundary vertex.  Then
$$\dim \ker \Delta^D=i(G,B).$$
In particular, if $i(G,B)=0$, then $\lambda_1>0$ which also implies that $\Delta^D$ is invertible.  Especially, when $G$ is connected, $\lambda_1>0$.
\end{prop}
\begin{proof}
Let $u\in \ker\Delta^D$. By \eqref{eq-int-by-part-Dirichlet}, we have
$$\vv<dE_0(u),dE_0(u)>=0.$$
So, $E_0(u)$ is constant on each connected component of $G$. Because $E_0(u)|_B=0$,  $E_0(u)=0$ on each connected component of $G$ that contains some boundary vertices. On those connected components of $G$ that contain no boundary vertex, the constant can be arbitrary. So,
$$\dim \ker \Delta^D=i(G,B).$$
\end{proof}

A real number $\nu$ is called a Neumann eigenvalue of $(G,B,m,w)$  if the following Neumann boundary value problem:
\begin{equation*}
\left\{\begin{array}{ll}\Delta u(x)=\nu u(x)&x\in\Omega\\
\frac{\p u}{\p n}(x)=0&x\in B
\end{array}\right.
\end{equation*}
has a nontrivial solution. The corresponding operator for Neumann  eigenvalues is the Neumann Laplacian operator introduced as follows.

For each $u\in \R^\Omega$, let $N_0(u)\in \R^V$ be its harmonic extension. That is
\begin{equation}\label{eq-N0}
\left\{\begin{array}{ll}\Delta N_0(u)(x)=0&x\in B\\
N_0(u)(x)=u(x)&x\in \Omega.
\end{array}\right.
\end{equation}
The Neumann eigenvalues of $(G,B,m,w)$ are the eigenvalues of the operator $\Delta^N:\R^\Omega\to\R^\Omega$ where
\begin{equation}\label{eq-Neu-Lap}
\Delta ^Nu=\Delta N_0(u)|_\Omega
\end{equation}
is called the Neumann Laplacian operator.

Note that harmonic extension is not unique without any connectivity assumption on $(G,B)$. However, the definition \eqref{eq-Neu-Lap} of $\Delta^N:\R^\Omega\to \R^\Omega$ is independent of the choices of the harmonic extension $N_0(u)$. In fact, consider all the connected components of $G$, on those connected components containing some interior vertices, the harmonic extension is uniquely determined by $u$ (see \cite{Gr}). So, $\Delta N_0(u)|_\Omega$ is uniquely determined by $u$. For simplicity, we also require that $N_0(u)=0$ on connected components of $G$ that contain no interior vertices. With this requirement on $N_0(u)$, $N_0:\R^\Omega\to \R^V$ is a well defined linear map.

Moreover, note that
\begin{equation}\label{eq-int-by-part-Neumann}
\vv<\Delta^N u,v>_\Omega=\vv<\Delta N_0(u),N_0(v)>=\vv<dN_0(u),dN_0(v)>
\end{equation}
by \eqref{eq-integration-by-part}. So $\Delta^N:\R^\Omega\to \R^\Omega$ is a nonnegative self-adjoint operator. Let
\begin{equation*}
0\leq \nu_1<\nu_2\leq\cdots\leq\nu_{|\Omega|}
\end{equation*}
be its eigenvalues. We also denote $\Delta^N$ as $\Delta^N_{G,B}$, and denote $\nu_i$ as $\nu_i(G,B)$ or $\nu_i(G,B,m,w)$ when it is necessary. When $i>|\Omega|$, we take the convention that $\nu_i=+\infty$.

By \eqref{eq-int-by-part-Neumann}, we have the following conclusion about the multiplicity of the zero eigenvalue for the Neumann Laplacian.
\begin{prop}\label{prop-ker-N}
Let $(G,B,m,w)$ be a weighted finite  graph with boundary. Let $j(G,B)$ be the number of connected components of $G$ that contain some interior vertices. Then,
$$\dim \ker \Delta^N=j(G,B).$$
In particular, if $G$ is connected with nonempty interior, then $\nu_2>0$.
\end{prop}
\begin{proof}
Let $u\in \ker \Delta ^N$. Then, by \eqref{eq-int-by-part-Neumann},
$$\vv<dN_0(u),dN_0(u)>=0.$$
So, $N_0(u)$ is constant on each connected components of $G$. Note that
$$N_0(u)|_\Omega=u$$
and the constants on each  connected components are arbitrary. So,
$$\dim\ker\Delta^N=j(G,B).$$
\end{proof}

A real number $\sigma$ is called a Steklov eigenvalue of $(G,B,m,w)$ if the following boundary value problem:
\begin{equation*}
\left\{\begin{array}{ll}\Delta u(x)=0&x\in\Omega\\
\frac{\p u}{\p n}(x)=\sigma u(x)&x\in B
\end{array}\right.
\end{equation*}
has a nonzero solution. The Steklov eigenvalues are the eigenvalues of the Steklov operator introduced as follows.

For any $u\in \R^B$, let $H_0(u)\in\R^V$ be the harmonic extension of $u$. That is,
\begin{equation*}
\left\{\begin{array}{ll}\Delta H_0(u)(x)=0&x\in\Omega\\
H_0(u)(x)=u(x)&x\in B.
\end{array}\right.
\end{equation*}
The Steklov eigenvalues are the eigenvalues of the Steklov operator $\Lambda:\R^B\to \R^B$ defined as
$$\Lambda(u)=\frac{\p H_0(u)}{\p n}.$$
The same as in the definition of Neumann operators, although harmonic extension is not unique in general, the Steklov operator is independent of the choices of harmonic extensions. Moreover, if we further require that $H_0(u)$ is zero on those connected components contain no boundary vertex, then $H_0:\R^B\to \R^V$ is a well defined linear map.

By \eqref{eq-integration-by-part}, it is clear that
$$\vv<\Lambda u,u>_B=\vv<\Delta H_0(u),H_0(u)>=\vv<dH_0(u),dH_0(u)>.$$
So, $\Lambda$ is a nonnegative self-adjoint operator, let
$$0=\sigma_1\leq\sigma_2\leq\cdots\leq \sigma_{|B|}.$$
be is eigenvalues.

We also denote $\Lambda$ as $\Lambda_{G,B}$ and denote $\sigma_i$ as
$\sigma_i(G,B)$ or $\sigma_i(G,B,m,w)$ if necessary. We take the convention that $\sigma_{i}=+\infty$ when $i>|B|$.

Note that in the general definition of graphs with boundary, one can interchange the roles of the interior and boundary. From this point of view, by the definition of Neumann eigenvalues and Steklov eigenvalues, one can see that they are dual of each other. More precisely,  one has the following relation between Steklov operators and Neumann operators:
$$\Lambda_{G,B}=\Delta^N_{G,\Omega},$$
and hence
$$\sigma_i(G,B)=\nu_i(G,\Omega)$$
for $i=1,2,\cdots,|B|$. So, any result on Neumann eigenvalues can be translated as a result on Steklov eigenvalues by interchanging the roles of $\Omega$ and $B$. For example, Proposition \ref{prop-ker-N} translated as a result on Steklov eigenvalues is as follows.
\begin{cor}\label{cor-ker-S}
Let $(G,B,m,w)$ be a weighted finite graph with boundary, let $b(G,B)$ be the number of connected components of $G$ containing some boundary vertices. Then
$$\dim\ker \Lambda=b(G,B).$$
In particular, if $G$ is connected, then $\sigma_2>0$.
\end{cor}
\begin{rem}
The conclusion of Corollary \ref{cor-ker-S} is slightly different with that of \cite[Proposition 3.2]{HHW}. This is because connected components that consist of only interior vertices were ignored in \cite[Proposition 3.2]{HHW}.
\end{rem}

Finally, recall the following theorem of Weyl (see Theorem 4.3.1 of \cite{HJ}) on the spectrum of the sum of two self-adjoint linear transformations which will be frequently used in Section 5.
\begin{thm}[Weyl]\label{thm-Weyl}
Let $V$ be an $n$ dimensional vector space equipped with an inner product. Let $A,B$ be two self-adjoint linear transformations of $V$. Then,
\begin{enumerate}
\item for $i=1,2,\cdots, n$ and $j=0,1,\cdots,n-i$,
$$\lambda_i(A+B)\leq \lambda_{i+j}(A)+\lambda_{n-j}(B)$$
with equality for some pair $i,j$ if and only if there is a nonzero vector $x\in V$ such that $(A+B)x=\lambda_i(A)x$, $Ax=\lambda_{i+j}(A+B)x$ and $Bx=\lambda_{n-j}(B)x$;
\item for $i=1,2,\cdots,n$ and $j=1,2,\cdots,i$,
$$\lambda_i(A+B)\geq \lambda_{i-j+1}(A)+\lambda_j(B)$$
with equality for some pair $i,j$ if and only if there is a nonzero vector $x\in V$ such that $(A+B)x=\lambda_i(A+B)x$, $Ax=\lambda_{i-j+1}(A)x$ and $Bx=\lambda_{j}(B)x$.
\end{enumerate}
Here $\lambda_i(A)$ means the $i^{\rm th}$ eigenvalue of $A$.
\end{thm}

\section{Relations of $\Delta^D,\Delta^N$ and $\Delta_\Omega$}
In this section, we derive some simple but useful relations among $\Delta_\Omega$, $\Delta^D$ and $\Delta^N$.
\begin{prop}\label{prop-D}
Let $(G,B,m,w)$ be a weighted finite graph with boundary. Then,
for any $u\in \R^\Omega$,
\begin{equation}\label{eq-Diri-Lap-2}
\Delta^D=\Delta_\Omega +\Deg_B\cdot I_\Omega
\end{equation}
where $I_\Omega:\R^\Omega\to \R^\Omega$ is the identity map.
\end{prop}
\begin{proof}For any $u\in \R^\Omega$ and $x\in\Omega$,
\begin{equation*}
\begin{split}
\Delta^Du(x)=&\Delta E_0(u)(x)\\
=&\frac{1}{m_x}\sum_{y\in V}(E_0(u)(x)-E_0(u)(y))w_{xy}\\
=&\frac{1}{m_x}\sum_{y\in \Omega}(u(x)-u(y))w_{xy}+\frac{1}{m_x}\sum_{y\in B}u(x)w_{xy}\\
=&\Delta_\Omega u(x)+\Deg_B(x)u(x).\\
\end{split}
\end{equation*}
This completes the proof of the proposition.
\end{proof}
Note that removing connected components of $G$ that contain no interior vertex will not influence the expression of $\Delta_\Omega$, $\Delta^D$ and $\Delta^N$ by definition. So, we may assume that any connected components of $G$ contain some interior vertex. By Proposition \ref{prop-ker-D}, this implies that $\lambda_1(G,\Omega)>0$ and that
$$L_B:=\Delta^D_{G,\Omega}:\R^B\to \R^B$$
is invertible. With the help of $L_B^{-1}$, we have the following expression of $N_0(u)$.
\begin{prop}\label{prop-N}
Let $(G,B,m,w)$ be a weighted finite graph with boundary such that every connected components of $G$ contains some interior vertices. Then, for any $u\in \R^\Omega$,
$$N_0(u)|_B=L_B^{-1}A_\Omega u.$$
In particular, when $E(B,B)=\emptyset$,
$$N_0(u)|_B=\Deg^{-1}\cdot A_\Omega u.$$
\end{prop}
\begin{proof}
For any $u\in \R^\Omega$ and $x\in B$, we have
\begin{equation*}
\begin{split}
0=&\Delta N_0(u)(x)\\
=&\frac{1}{m_x}\sum_{y\in V}\left(N_0(u)(x)-N_0(u)(y)\right)w_{xy}\\
=&\frac{1}{m_x}\sum_{y\in B}\left(N_0(u)(x)-N_0(u)(y)\right)w_{xy}+\frac{1}{m_x}\sum_{y\in \Omega}\left(N_0(u)(x)-u(y)\right)w_{xy}\\
=&\Delta_BN_0(u)(x)+\Deg_\Omega(x)N_0(u)(x)-A_\Omega u(x).
\end{split}
\end{equation*}
By \eqref{eq-Diri-Lap-2}, we have
$$L_BN_0(u)=\Delta_BN_0(u)+\Deg_\Omega N_0(u).$$
So,
$$N_0(u)|_B=L_B^{-1}A_\Omega u.$$

Moreover, when $E(B,B)=\emptyset$, $\Delta_B=0$. So, we have
$$N_0(u)(x)=\frac{1}{\Deg_\Omega(x)}A_\Omega u(x)=\frac{1}{\Deg(x)}A_\Omega u(x),\ \forall x\in B$$
by that $E(B,B)=\emptyset$. This completes the proof of the proposition.
\end{proof}
By the help of the expression of $N_0(u)$ in Proposition \ref{prop-N}, we have the following relations among $\Delta^N$, $\Delta_\Omega$ and $\Delta^D$.
\begin{prop}\label{prop-N-Lap}
Let $(G,B,m,w)$ be a weighted finite graph with boundary such that every connected component of $G$ contains some interior vertices. Then,
\begin{equation}\label{eq-N-L-D-1}
\begin{split}
\Delta^N=&\Delta_\Omega +\Deg_B\cdot I_\Omega-A_BL_B^{-1}A_\Omega \\
=&\Delta^D -A_BL_B^{-1}A_\Omega .
\end{split}
\end{equation}
In particular, when $E(B,B)=\emptyset$,
\begin{equation}\label{eq-N-L-D-2}
\begin{split}
\Delta^N=&\Delta_\Omega +\Deg_B\cdot u-A_B\Deg^{-1}A_\Omega\\
=&\Delta^D -A_B\Deg^{-1}A_\Omega .
\end{split}
\end{equation}
\end{prop}
\begin{proof}
For any $u\in\R^\Omega$ and $x\in\Omega$, by Proposition \ref{prop-D} and Proposition \ref{prop-N},
\begin{equation*}
\begin{split}
\Delta^N u(x)=&\frac{1}{m_x}\sum_{y\in V}( N_0(u)(x)- N_0(u)(y))w_{xy}\\
=&\frac{1}{m_x}\sum_{y\in \Omega}(u(x)-u(y))w_{xy}+\frac{1}{m_x}\sum_{y\in B}(u(x)-N_0(u)(y))w_{xy}\\
=&\Delta_\Omega u+\Deg_B(x)u(x)-\frac{1}{m_x}\sum_{y\in B}L_B^{-1}A_\Omega u(y)w_{xy}\\
=&\Delta_\Omega u+\Deg_B(x)u(x)-(A_BL_B^{-1} A_\Omega) u(x)\\
=&\left(\Delta^D-A_BL_B^{-1}A_\Omega\right)u(x).
\end{split}
\end{equation*}
This completes the proof of the proposition.
\end{proof}
Translating Proposition \ref{prop-N-Lap} into Steklov operators, we have the following conclusion.
\begin{cor}\label{cor-S}
Let $(G,B,m,w)$ be a weighted finite graph with boundary such that every connected component of $G$ contains some boundary vertices. Then,
\begin{equation}\label{eq-S-L-1}
\begin{split}
\Lambda=&\Delta_B +\Deg_\Omega\cdot I_B-A_\Omega L_\Omega^{-1}A_B \\
=&L_B -A_\Omega L_\Omega^{-1}A_B.
\end{split}
\end{equation}
In particular, when $E(B,B)=\emptyset$,
\begin{equation}\label{eq-S-L-2}
\begin{split}
\Lambda=\Deg_\Omega\cdot I_B-A_\Omega L_\Omega^{-1}A_B.\\
\end{split}
\end{equation}
Here $L_\Omega=\Delta^D$.
\end{cor}
\begin{rem}
The identity \eqref{eq-S-L-2} is essentially Corollary 2.1 of \cite{HHW}. So, the identity \eqref{eq-S-L-1} can be viewed as an extension of \cite[Corollary 2.1]{HHW} for more general graphs with boundary.
\end{rem}
\section{Comparison of the Neumann eigenvalues and Laplacian eigenvalues and its rigidity}
In this section, we prove the the comparison of Neumann eigenvalues and Laplacian eigenvalues and its rigidity.

We first prove Theorem \ref{thm-Neumann-Laplacian}, a comparison of Neumann  eigenvalues and Laplacian eigenvalues.
\begin{proof}[Proof of Theorem \ref{thm-Neumann-Laplacian}] Let $v_1=1,v_2,\cdots, v_{|\Omega|}\in \R^\Omega$ be eigenfunctions of $\nu_1=0,\nu_2,\cdots,\nu_{|\Omega|}$ respectively such that $$\vv<v_i,v_j>_\Omega=0$$ for $i\neq j$. Let $u_1=1,u_2,\cdots, u_{|V|}\in \R^V$ be eigenfunctions of $\mu_1,\mu_2,\cdots,\mu_{|V|}$ respectively such that $$\vv<u_i,u_j>=0$$
for $i\neq j$. For each $i\geq 2$, let $v=c_1v_1+c_2v_2+\cdots+c_iv_i$ with $c_1,c_2,\cdots,c_i$ not all zero, such that
\begin{equation}\label{eq-linear-res}
\vv<N_0(v), u_j>=0
\end{equation}
for $j=1,2,\cdots, i-1$. This can be done because \eqref{eq-linear-res} for $j=1,2,\cdots, i-1$ form a homogeneous linear system with $i-1$ equations and $i$ unknowns which will certainly have nonzero solutions. Then,
\begin{equation}\label{eq-nu_i-mu_i}
\nu_i\geq\frac{\vv<dN_0(v),dN_0(v)>}{\vv<v,v>_\Omega}\geq\frac{\vv<dN_0(v),dN_0(v)>}{\vv<N_0(v),N_0(v)>}\geq \mu_i.
\end{equation}
If $\nu_i=\mu_i>0$, then all the inequalities become equalities in \eqref{eq-nu_i-mu_i}. Hence, $N_0(v)$ is simultaneously a Neumann Laplacian eigenfunction and a Laplacian eigenfunction for $\nu_i=\mu_i$. Moreover, since $\mu_i>0$,
$$\vv<dN_0(v),dN_0(v)>>0.$$
Thus, $N_0(v)|_B=0$. Note that $\frac{\p N_0(v)}{\p n}=0$ by definition of $N_0$. This completes the proof of the theorem by letting $u=N_0(v)$.
\end{proof}
We next prove Theorem \ref{thm-Neu-Lap-rigidity}, a general rigidity for \eqref{eq-Neumann-Laplacian}.

\begin{proof}[Proof of Theorem \ref{thm-Neu-Lap-rigidity}]
If the equality of \eqref{eq-Neumann-Laplacian} holds for all $i=1,2,\cdots, |\Omega|$, we claim that \emph{there are $\tilde u_1=1, \tilde u_2,\cdots, \tilde u_{|\Omega|}\in \R^V$ such that
\begin{enumerate}
\item $\Delta \tilde u_i=\mu_i \tilde u_i$ for $i=1,2,\cdots,|\Omega|$;
\item $\Delta^N(\tilde u_i|_\Omega)=\nu_i\tilde u_i|_\Omega$ for $i=1,2,\cdots,|\Omega|$;
\item $\tilde u_i|_B=\frac{\p\tilde u_i}{\p n}=0$ for $i=2,3,\cdots, |\Omega|$;
\item $\vv<\tilde u_i,\tilde u_j>_\Omega=\vv<\tilde u_i,\tilde u_j>=0$ for $1\leq j<i\leq |\Omega|$.
\end{enumerate}}
In fact, suppose that $\tilde u_1=1,\tilde u_2,\cdots,\tilde u_{i-1}$ have been constructed for some $i\geq 2$. Let $v=c_1v_1+c_2v_2+\cdots+c_iv_i$ with $c_1,c_2,\cdots,c_i$ not all zero such that
$$\vv<N_0(v),\tilde u_j>=0$$
for $j=1,2,\cdots,i-1$. Here $v_1,v_2,\cdots,v_{|\Omega|}$ are the same as in the proof of Theorem \ref{thm-Neumann-Laplacian} Then, using the same argument as in \eqref{eq-nu_i-mu_i}, we know that $\tilde u_i:=N_0(v)$ will satisfy the requirements (1)-(4) above, since $\nu_i=\mu_i>0$ for $i\geq 2$.

Then, for any $v\in \R^\Omega$ with
\begin{equation}\label{eq-N-L-1}
0=\vv<v,1>_\Omega=\sum_{y\in \Omega}v(y)m_y.
\end{equation}
We know that
\begin{equation*}
v=c_2\tilde u_2|_\Omega+\cdots c_{|\Omega|}\tilde u_{|\Omega|}|_{\Omega}
\end{equation*}
for some $c_2,c_2,\cdots,c_{|\Omega|}\in\R$. Thus, by (3) of the claim,
\begin{equation*}
N_0(v)=c_2\tilde u_2+\cdots+c_{|\Omega|}\tilde u_{\Omega}
\end{equation*}
which implies that $N_0(v)|_B=0$ by (3) of the claim again. By Proposition \ref{prop-N},
\begin{equation*}
\frac1{m_x}\sum_{y\in\Omega}v(y)w_{xy}=A_\Omega v(x)=\left(L_B N_0(v)|_B\right)(x)=0
\end{equation*}
for any $x\in B$. By comparing this to \eqref{eq-N-L-1}, there is a nonnegative function $\rho\in \R^B$ such that
$$w_{xy}=\rho_xm_xm_y$$ for any $x\in B$ and $y\in\Omega$.

Conversely, if $w_{xy}=\rho_xm_xm_y$ for any $x\in B$ and $y\in\Omega$,
then for any $v\in \R^\Omega$  with $\vv<v,1>_\Omega=0$, we have
$$A_\Omega v(x)=\rho_x\vv<v,1>_\Omega=0,\ \forall x\in B.$$
So, by Proposition \ref{prop-N},
$$N_0(v)|_B=L_B^{-1}A_\Omega v\equiv 0.$$
If $v$ is further an eigenfunction of $\Delta^N$ with eigenvalue $\nu_i$, then by that
$$\Delta N_0(v)|_\Omega=\Delta^N v=\nu_i v$$
and
$$\Delta N_0(v)|_B=0=\nu_i N_0(v)|_B,$$
we know that $N_0(v)$ is also an eigenfunction of $\Delta$ with the same eigenvalue $\nu_i$.

Moreover, for any $u\in \R^\Omega$ with $\vv<u,1>_\Omega=0$, since $N_0(u)|_B=0$, by Proposition \ref{prop-D},
\begin{equation*}
\Delta^N u=\Delta N_0(u)=\Delta E_0(u)|_\Omega=\Delta^D u=\Delta_\Omega u+\Deg_b \cdot u.
\end{equation*}
So,
\begin{equation*}
\mu_{|\Omega|}=\nu_{|\Omega|}=\mu_{|\Omega|}(\Omega)+\Deg_b.
\end{equation*}
Hence, the equality $\nu_i=\mu_i$ for $i=1,2,\cdots, |\Omega|$ holds if and only if
\begin{equation}\label{eq-nu-mu}
\frac{\vv<du,du>}{{\vv<u,u>}}\geq \mu_{|\Omega|}(\Omega)+\Deg_b
\end{equation}
for any nonzero $u\in\R^V$ with $\vv<u,\tilde u_i>=0$ for $i=1,2,\cdots, |\Omega|$. Because
$$\vv<u,\tilde u_i>_\Omega=\vv<u,\tilde u_i>=0$$
for $i=2,3,\cdots, |\Omega|$, we know that $u|_\Omega$ must be a constant $c$.

When $c=0$, let $f=u|_B$, then $\vv<f,1>_B=\vv<u,1>=0$ and by \eqref{eq-nu-mu},
\begin{equation*}
\frac{\vv<df,df>_B+V_\Omega\vv<\rho f,f>_B}{\vv<f,f>_B}=\frac{\vv<du,du>}{\vv<u,u>}\geq \mu_{|\Omega|}(\Omega)+\Deg_b.
\end{equation*}
That is,
\begin{equation}\label{eq-f-1}
\vv<df,df>_B+V_\Omega\vv<\rho  f,f>_B-\left(\mu_{|\Omega|}(\Omega)+\Deg_b\right)\vv<f,f>_B\geq 0
\end{equation}
for any $f\in \R^B$ with $\vv<f,1>_B=0$.

When $c\neq 0$, we can assume that $c=1$. Let $f=u|_B+\frac{V_\Omega}{V_B}$. Then, by that $\vv<u,1>=0$, we have $\vv<f,1>_B=0$. Moreover, by \eqref{eq-nu-mu},
\begin{equation*}
\begin{split}
&\frac{\vv<df,df>_B+V_\Omega\left(\vv<\rho f,f>_B-\frac{2V_G}{V_B}\vv<\rho, f>_B+\frac{V_G^2}{V_B^2}\Deg_b\right)}{\vv<f,f>_B+\frac{V_\Omega V_G}{V_B}}\\
=&\frac{\vv<df,df>_B+V_\Omega\vv<\rho(f-\frac{V_G}{V_B}),f-\frac{V_G}{V_B}>_B}{V_\Omega+\vv<f-\frac{V_\Omega}{V_B},f-\frac{V_\Omega}{V_B}>_B}\\
=&\frac{\vv<du,du>}{\vv<u,u>}\\
\geq&\mu_{|\Omega|}(\Omega)+\Deg_b.
\end{split}
\end{equation*}
That is,
\begin{equation}\label{eq-Neu-Lap-3}
\begin{split}
\vv<df,df>_B+&V_\Omega\vv<\rho f,f>_B-\left(\mu_{|\Omega|}(\Omega)+\Deg_b\right)\vv<f,f>_B\\
&-2\frac{V_\Omega V_G}{V_B}\vv<\rho, f>_B+\frac{V_\Omega V_G}{V_B}\left(\frac{V_\Omega}{V_B}\Deg_b-\mu_{|\Omega|}(\Omega)\right)\geq0
\end{split}
\end{equation}
for any $f\in \R^B$ with $\vv<f,1>_B=0$. Let $f=0$ in \eqref{eq-Neu-Lap-3}. We get
\begin{equation}\label{eq-f-2}
\mu_{|\Omega|}(\Omega)\leq \frac{V_\Omega}{V_B}\Deg_b.
\end{equation}
 Moreover, replacing $f$ by $\lambda f$ in \eqref{eq-Neu-Lap-3}, we have
\begin{equation*}
\begin{split}
&\left(\vv<df,df>_B+V_\Omega\vv<\rho f,f>_B-\left(\mu_{|\Omega|}(\Omega)+\Deg_b\right)\vv<f,f>_B\right)\lambda^2\\
&\ \ \ \ \ \ \ \ \ \ \ \ -2\frac{V_\Omega V_G}{V_B}\vv<\rho, f>_B\lambda+\frac{V_\Omega V_G}{V_B}\left(\frac{V_\Omega}{V_B}\Deg_b-\mu_{|\Omega|}(\Omega)\right)\geq0
\end{split}
\end{equation*}
for any $\lambda\in \R$. Then,
\begin{equation*}
\begin{split}
&\left(2\frac{V_\Omega V_G}{V_B}\vv<\rho, f>_B\right)^2\\
\leq& 4\frac{V_\Omega V_G}{V_B}\left(\vv<df,df>_B+V_\Omega\vv<\rho f,f>_B-\left(\mu_{|\Omega|}(\Omega)+\Deg_b\right)\vv<f,f>_B\right)\left(\frac{V_\Omega}{V_B}\Deg_b-\mu_{|\Omega|}(\Omega)\right)
\end{split}
\end{equation*}
which is equivalent to
\begin{equation}\label{eq-f-3}
\begin{split}
&\left(\frac{V_\Omega}{V_B}\Deg_b-\mu_{|\Omega|}(\Omega)\right)\left(\vv<df,df>_B+V_\Omega\vv<\rho f,f>_B-\left(\mu_{|\Omega|}(\Omega)+\Deg_b\right)\vv<f,f>_B\right)\\
&-\frac{V_\Omega V_G}{V_B}\vv<\rho,f>_B^2\geq 0
\end{split}
\end{equation}
for any $f\in \R^B$ with $\vv<f,1>_B=0$.

Conversely, it is not hard to see that the combination of \eqref{eq-f-1}, \eqref{eq-f-2} and \eqref{eq-f-3} implies \eqref{eq-nu-mu}.

When $\rho$ is constant, substituting this into \eqref{eq-f-1}, \eqref{eq-f-2} and \eqref{eq-f-3}, we find that they are equivalent to that
\begin{equation}\label{eq-df}
\vv<df,df>_B+\left(\rho(V_\Omega-V_B)-\mu_{|\Omega|}(\Omega)\right)\vv<f,f>_B\geq 0
\end{equation}
for all $f\in \R^B$ with $\vv<f,1>_B=0$,
and
\begin{equation*}
\mu_{|\Omega|}(\Omega)\leq \rho V_\Omega
\end{equation*}
by noting that $\Deg_b=\rho V_B$ in this case. Note that \eqref{eq-df} is automatical true if $|B|=1$ and when $|B|\geq 2$, \eqref{eq-df} is equivalent to that
$$\mu_{|\Omega|}(\Omega)\leq \mu_2(B)+\rho(V_\Omega-V_B).$$
So, when $\rho$ is constant, the combination of \eqref{eq-f-1}, \eqref{eq-f-2} and \eqref{eq-f-3} is equivalent to
$$\mu_{|\Omega|}(\Omega)\leq \min\left\{\rho V_\Omega,\mu_2(B)+\rho(V_\Omega-V_B)\right\}.$$

When $\rho$ is nonconstant, let $f\in \R^B$  be such that $\vv<1,f>_B=0$ and $\vv<\rho,f>_B\neq0$. Substituting $f$ into \eqref{eq-f-3}, we get
$$\mu_{|\Omega|}(\Omega)<\frac{V_\Omega}{V_B}\Deg_b.$$
Moreover, \eqref{eq-f-3} can be rewritten as
\begin{equation}\label{eq-r-2}
\begin{split}
&\vv<df,df>_B+V_\Omega\vv<\rho f,f>_B-\left(\mu_{|\Omega|}(\Omega)+\Deg_b\right)\vv<f,f>_B\\
\geq& \frac{V_\Omega V_G}{V_\Omega\Deg_b-V_B\mu_{|\Omega|}(\Omega)}\vv<\rho,f>_B^2.
\end{split}
\end{equation}
which is stronger than \eqref{eq-f-1}. So, we only need to require \eqref{eq-r-2} in this case.

This completes the proof of the theorem.
\end{proof}
We next come to prove Corollary \ref{cor-rigidity-unit-weight} and Corollary \ref{cor-rigidity-normalized-weight}, the rigidity of \eqref{eq-Neumann-Laplacian} when the graph is equipped with the unit weight or normalized weight.
\begin{proof}[Proof of Corollary \ref{cor-rigidity-unit-weight} ]
If the equalities of \eqref{eq-Neumann-Laplacian} holds for $i=1,2,\cdots,|\Omega|$, by Theorem \ref{thm-Neu-Lap-rigidity}, there is a nonnegative function $\rho\in \R^B$, such that
$$w_{xy}=\rho_xm_xm_y$$
for all $x\in B$ and $y\in \Omega$. Because $x\in B$ is adjacent to some $y\in \Omega$ and the graph is equipped with unit weight, we have
$$\rho\equiv 1.$$
So, any boundary vertex is adjacent to any interior vertex. Substituting $\rho\equiv1$ into \eqref{eq-com-Neu-Lap-rho-const}, we get
\begin{equation*}
  \mu_{|\Omega|}(\Omega)\leq\min\{|\Omega|, \mu_2(B)+|\Omega|-|B|\}=\mu_2(B)+|\Omega|-|B|
\end{equation*}
when $|B|\geq 2$, by noting that $\mu_2(B)\leq |B|$ (see Fiedler \cite{FI}). When $|B|=1$, $\mu_2(B)=+\infty$, by \eqref{eq-com-Neu-Lap-rho-const},
$$\mu_{|\Omega|}(\Omega)\leq |\Omega|.$$

Conversely, because every boundary vertex is adjacent to every interior vertex, we have
$$w_{xy}=m_xm_y$$
for any $x\in B$ and $y\in \Omega$. So $\rho\equiv1$. Then, note that (2) of Corollary \ref{cor-rigidity-unit-weight} is equivalent to \eqref{eq-com-Neu-Lap-rho-const}. So, by Theorem \ref{thm-Neu-Lap-rigidity}, the equalities of \eqref{eq-com-Neu-Lap-rho-const} hold for $i=1,2,\cdots,|\Omega|$.

Moreover, recall the following estimate (see \cite{FI} for example):
$$\mu_{|\Omega|}(\Omega)\leq 2\Delta(G|_\Omega).$$
So, if $$2\Delta(G|_\Omega)\leq |\Omega|-|B|,$$
then (2) of Corollary \ref{cor-rigidity-unit-weight} holds. This completes the proof of Corollary \ref{cor-rigidity-unit-weight}.
\end{proof}
\begin{proof}[Proof of  Corollary \ref{cor-rigidity-normalized-weight}]
Suppose the equalities of \eqref{eq-Neumann-Laplacian} hold for $i=1,2,\cdots,|\Omega|$. By Theorem \ref{thm-Neu-Lap-rigidity}, for any $x\in B$,
\begin{equation*}
1=\Deg(x)=\frac{1}{m_x}\sum_{y\in \Omega}w_{xy}=\rho_x V_\Omega,
\end{equation*}
since $E(B,B)=\emptyset$. So, for any $x\in B$, $\rho_x=\frac{1}{V_\Omega}$. Substituting this into \eqref{eq-com-Neu-Lap-rho-const}, we get
\begin{equation*}
\mu_{|\Omega|}(\Omega)\leq \frac{V_\Omega-V_B}{V_\Omega}
\end{equation*}
when $|B|\geq 2$, since $\mu_2(B)=0$ by $E(B,B)=\emptyset$. When $|B|=1$, since $\mu_2(B)=+\infty$, by \eqref{eq-com-Neu-Lap-rho-const},
$$\mu_{|\Omega|}(\Omega)\leq 1.$$

The converse is clearly true by Theorem \ref{thm-Neu-Lap-rigidity}. This completes the proof of Corollary \ref{cor-rigidity-normalized-weight}.
\end{proof}

\section{Comparisons of the eigenvalues for $\Delta^N$,$\Delta^D$ and $\Delta_\Omega$}
In this section, we obtain the comparisons of Dirichlet eigenvalues, Neumann eigenvalues and Laplacian eigenvalues of the interior subgraph, and their rigidities.

We first prove Theorem \ref{thm-D-L}, a comparison of the Dirichlet eigenvalues and the Laplacian eigenvalues of the interior subgraph.
\begin{proof}[Proof of Theorem \ref{thm-D-L}]
By Proposition \ref{prop-D}, we have
$$\Delta^D =\Delta_\Omega +\Deg_B\cdot I_\Omega$$
Then, the conclusions come from Theorem \ref{thm-Weyl} directly by noting that
$$\lambda_i(\Deg_B\cdot I_\Omega)=\Deg_B^{(i)}$$
for $i=1,2,\cdots,|\Omega|$.
\end{proof}

We next prove Theorem \ref{thm-N-L}, a comparison of Neumann eigenvalues and Laplacian eigenvalues of the interior subgraph.
\begin{proof}[Proof of Theorem \ref{thm-N-L}]
For any $u\in \R^\Omega$, by \eqref{eq-int-by-part-Neumann}, we have
\begin{equation*}
\vv<\Delta^N u,u>_\Omega=\vv<dN_0(u),dN_0(u)>\geq \vv<du,du>_\Omega=\vv<\Delta_\Omega u,u>_\Omega.
\end{equation*}
So, for any $i=1,2,\cdots,|\Omega|$, $\nu_i\geq \mu_i(\Omega)$.

Moreover,  $\nu_i=\mu_i(\Omega)$ for all $i=1,2,\cdots,|\Omega|$ if and only if  for any $u\in \R^\Omega$,
$$\vv<dN_0(u),dN_0(u)>=\vv<du,du>_\Omega.$$
This is equivalent to
$$\vv<dN_0(u),dN_0(u)>_B=\vv<dN_0(u),dN_0(u)>_{E(B,\Omega)}=0$$
and that $N_0(u)$ is constant on each connected component of $N_{G}(B)$ for any $u\in \R^\Omega$. Because $u\in \R^\Omega$ is arbitrary, this is equivalent to that each connected component of $N_G(B)$ contains at most one interior vertex.
\end{proof}
By Proposition \ref{prop-N-Lap}, the identity on the Neumann Laplacian operator and the Laplacian operator on the interior subgraph, we have the following corollary.
\begin{cor}\label{cor-Deg}
Let $(G,B,m,w)$ be a weighted finite graph with boundary such that each connected component contains some interior vertices. Then,
\begin{enumerate}
\item $\Deg_B\cdot I_\Omega-A_BL_B^{-1}A_\Omega$ is nonnegative;
\item for $i=1,2,\cdots,|\Omega|$,
$$\lambda_i(A_BL_B^{-1}A_\Omega)\leq \Deg_B^{(i)}$$
 with equalities for all $i=1,2,\cdots,|\Omega|$ if and only if each connected component of $N_G(B)$ contains at most one interior vertex.
\end{enumerate}
\end{cor}
\begin{proof}
By Proposition \ref{prop-N-Lap}, \eqref{eq-integration-by-part} and \eqref{eq-int-by-part-Neumann}, for any $u\in \R^\Omega$
\begin{equation*}
\begin{split}
\vv<\Deg_B\cdot u-A_BL_B^{-1}A_\Omega u,u>_\Omega=&\vv<\Delta^N u-\Delta_\Omega u,u>_\Omega\\
=&\vv<dN_0(u),N_0(u)>-\vv<du,du>_\Omega\geq0.
\end{split}
\end{equation*}
So $\Deg_B\cdot I_\Omega-A_BL_B^{-1}A_\Omega$ is nonnegative. This proves (1).

The inequalities in (2) comes directly from (1). When the equalities in (2) all hold, we know that $\Delta^N=\Delta_\Omega$. So, each connected components of $N_G(B)$ contains at most one interior vertex. This completes the proof of the corollary.
\end{proof}
Translating Theorem \ref{thm-N-L} to Steklov eigenvalues, we have the following conclusion.
\begin{cor}\label{cor-Steklov-0}
Let $(G,B,m,w)$ be a weighted finite graph with boundary. Then,
$$\sigma_i(G,B)\geq \mu_i(B)$$
for $i=1,2,\cdots,|B|$, with equalities for $i=1,2,\cdots,|B|$ if and only if each connected components of $N_G(\Omega)$ contains at most one boundary vertex.
\end{cor}

We next prove Theorem \ref{thm-N-D}, a comparison of Neumann eigenvalues and Dirichlet eigenvalues.
\begin{proof}[Proof of Theorem \ref{thm-N-D}]
For any $u\in \R^\Omega$, by the Dirichlet Principle (see \cite{Gr}), we have
$$\vv<dE_0(u),dE_0(u)>\geq \vv<dN_0(u),dN_0(u)>.$$
So, by \eqref{eq-int-by-part-Dirichlet} and \eqref{eq-int-by-part-Neumann}, we have
$$\vv<\Delta^D u,u>_\Omega=\vv<dE_0(u),dE_0(u)>\geq \vv<dN_0(u),dN_0(u)>=\vv<\Delta^N u,u>_\Omega.$$
Thus, $\lambda_i\geq \nu_i$ for $i=1,2,\cdots,|\Omega|$.
\end{proof}
Translating Theorem \ref{thm-N-D} into Steklov eigenvalues will give us the following estimate on Steklov eigenvalues which can be viewed as an extension of the upper bound on Dirichlet-to-Neumann map in \cite[Proposition 3.1]{HHW}.
\begin{cor}\label{cor-Steklov-1}
Let $(G,B,m,w)$ be a weighted finite graph. Then,
$$\sigma_i(G,B)\leq \mu_{i+j}(B)+\Deg_\Omega^{(|B|-j)}$$
for any $i=1,2,\cdots, |B|$ and $j=0,1,\cdots,|B|-i$. In particular, if $E(B,B)=\emptyset$, then
$$\sigma_i(G,B)\leq \Deg_\Omega^{(i)}.$$
Here
$$\Deg_\Omega^{(1)}\leq \Deg_\Omega^{(2)}\leq\cdots\leq \Deg_\Omega^{(|B|)}$$
is the rearrangement of $\Deg_\Omega(x)$ with $x$ going through $B$.
\end{cor}
\begin{proof}
By Theorem \ref{thm-N-D} and (1) of Theorem \ref{thm-D-L}, we have
\begin{equation*}
\begin{split}
\sigma_i(G,B)=\nu_i(G,\Omega)\leq\lambda_i(G,\Omega)\leq \mu_{i+j}(B)+\Deg_\Omega^{(|B|-j)}
\end{split}
\end{equation*}
for any $i=1,2,\cdots, |B|$ and $j=0,1,\cdots,|B|-i$.

When $E(B,B)=\emptyset$, we have $\mu_{|B|}(B)=0$. So, by letting $j=|B|-i$ in the last inequality, we get the inequality $\sigma_i(G,B)\leq \Deg_\Omega^{(i)}$ for $i=1,2,\cdots,|B|$.
\end{proof}
We next prove Theorem \ref{thm-N-D-sharp} and Theorem \ref{thm-N-D-sharp-1}, sharp comparisons of Neumann eigenvalues and Dirichlet eigenvalues.
\begin{proof}[Proof of Theorem \ref{thm-N-D-sharp}] By Proposition \ref{prop-N-Lap},
$$\Delta^D=\Delta^N+A_BL_B^{-1}A_\Omega.$$
Then, the conclusion comes directly comes from this identity.
\end{proof}
\begin{proof}[Proof of Theorem \ref{thm-N-D-sharp-1}]
The inequalities comes directly from Theorem \ref{thm-N-D-sharp} and Proposition \ref{prop-N}. We only need to show the rigidity. The equalities of \eqref{eq-com-Neu-Diri} hold for $i=1,2,\cdots,|\Omega|$ if and only if  $A_B\Deg^{-1} A_\Omega$ is a scalar operator. By direct computation, this is equivalent to the follows:
 \begin{equation}\label{eq-Neu-Diri-1}
 \sum_{x\in B}\frac{w_{xz}^2}{m_z\sum_{y\in\Omega}w_{xy}}:=s
 \end{equation}
 is independent of $z\in \Omega$,
 and
 \begin{equation}\label{eq-Neu-Diri-2}
 \sum_{x\in B}\frac{w_{xz}w_{xy}}{m_z\sum_{\xi\in\Omega}w_{x\xi}}=0
 \end{equation}
 for any $y\neq z\in\Omega$. Because $E(B,B)=\emptyset$ and $G$ is connected, each boundary vertex must be adjacent to some interior vertex. So, $s>0$. Then, by \eqref{eq-Neu-Diri-1}, we know that each interior vertex must be adjacent to some boundary vertex. Finally, it is clear that \eqref{eq-Neu-Diri-2} is equivalent to that each boundary vertex is adjacent to only one interior vertex. So, the there is a surjective map $\varphi:B\to \Omega$ such that $x\sim \varphi(x)$. Then,
 by \eqref{eq-Neu-Diri-1}, for any $z\in\Omega$,
 \begin{equation*}
 \begin{split}
 s=\sum_{x\in B}\frac{w_{xz}^2}{m_z\sum_{y\in\Omega}w_{xy}}
 =\sum_{x\in \varphi^{-1}(z)}\frac{w_{xz}^2}{m_z\sum_{y\in\Omega}w_{xy}}
 =\sum_{x\in \varphi^{-1}(z)}\frac{w_{xz}^2}{m_zw_{xz}}
 =\Deg_B(z).
 \end{split}
 \end{equation*}
 So, $\Deg_B(z)$ is independent of $z\in\Omega$.

 Conversely, it is not hard to show \eqref{eq-Neu-Diri-1} and \eqref{eq-Neu-Diri-2} assuming (1) and (2) of Theorem \ref{thm-N-D-sharp-1}. This completes the proof of the theorem.
\end{proof}

In fact, by using Theorem \ref{thm-Weyl} similarly as in the proof of Theorem \ref{thm-D-L}, we can get the following comparison of Neumann eigenvalues and Dirichlet eigenvalues which is an extension of Theorem \ref{thm-N-D-sharp-1}.
\begin{thm}\label{thm-N-D-sharp-2}
Let $(G,B,m,w)$ be a  weighted finite graph with boundary such that each connected component of $G$ contains some interior vertices. Then
\begin{enumerate}
\item for $i=1,2,\cdots,|\Omega|$, $j=0,1,\cdots,|\Omega|-i$,
$$\lambda_{i}\leq \nu_{i+j}+\lambda_{|\Omega|-j}(A_BL_B^{-1}A_\Omega)$$
with equality for some pair $i,j$ if and only if there is a nonzero function $u\in \R^\Omega$ such that $\Delta^Du=\lambda_i u$, $\Delta^Nu=\nu_{i+j} u$ and $A_BL_B^{-1}A_\Omega u=\lambda_{|\Omega|-j}(A_BL_B^{-1}A_\Omega)u$;
\item for $i=1,2,\cdots,|\Omega|$, $j=1,2,\cdots,i$,
$$\lambda_i\geq \nu_{i-j+1}+\lambda_j(A_BL_B^{-1}A_\Omega)$$
with equality for some pair $i,j$ if and only if there is a nonzero function $u\in \R^\Omega$ such that $\Delta^Du=\lambda_i u$, $\Delta^Nu=\nu_{i-j+1} u$ and $A_BL_B^{-1}A_\Omega u=\lambda_{j}(A_BL_B^{-1}A_\Omega)u$.
\end{enumerate}
\end{thm}
\begin{proof}
By Proposition \ref{prop-N-Lap}, we have
$$\Delta^D=\Delta^N+A_BL_B^{-1}A_\Omega.$$
Then, the conclusions follow directly by using Theorem \ref{thm-Weyl}.
\end{proof}

Translating the last result into Steklov eigenvalues, we get the following  estimate on Steklov eigenvalues.
\begin{cor}\label{cor-Steklov-2}
Let $(G,B,m,w)$ be a  weighted finite graph with boundary such that each connected component of $G$ contains some boundary vertices. Then
\begin{equation}\label{eq-S-1}
\lambda_{i-j}(G,\Omega)-\lambda_{|B|-j}(A_\Omega L_\Omega^{-1}A_B)\leq \sigma_i(G,B)\leq \lambda_{i+k-1}(G,\Omega)-\lambda_k(A_\Omega L_\Omega^{-1}A_B)
\end{equation}
for $i=1,2,\cdots,|B|$, $j=0,1,\cdots,i-1$ and $k=1,2,\cdots,|B|-i+1$. In particular, when $E(B,B)=\emptyset$, we have
\begin{equation}\label{eq-S-2}
\Deg_\Omega^{(i-j)}-\lambda_{|B|-j}(A_\Omega L_\Omega^{-1}A_B)\leq \sigma_i(G,B)\leq \Deg_\Omega^{(i+k-1)}-\lambda_k(A_\Omega L_\Omega^{-1}A_B)
\end{equation}
for $i=1,2,\cdots,|B|$, $j=0,1,\cdots,i-1$ and $k=1,2,\cdots,|B|-i+1$. Here $L_\Omega=\Delta^D$.
\end{cor}
\begin{proof}
By (1) of Theorem \ref{thm-N-D-sharp-2}, we have
\begin{equation*}
\lambda_{i-j}(G,\Omega)\leq \nu_{i}(G,\Omega)+\lambda_{|B|-j}(A_\Omega L_{\Omega}^{-1}A_B)=\sigma_i(G,B)+\lambda_{|B|-j}(A_\Omega L_{\Omega}^{-1}A_B).
\end{equation*}
for $i=1,2,\cdots,|B|$, $j=0,1,\cdots,i-1$. This gives us the lower bound in \eqref{eq-S-1}. The upper bound of \ref{eq-S-1} can be proved similarly by using (2) of Theorem \ref{thm-N-D-sharp-2}.

When $E(B,B)=\emptyset$, by Proposition \ref{prop-D}, we have
$$\Delta_{G,\Omega}^Du=\Delta_Bu+\Deg_\Omega\cdot u=\Deg_\Omega\cdot u,\ \forall u\in \R^B.$$
So, $\lambda_i(G,\Omega)=\Deg_\Omega^{(i)}$ for $i=1,2,\cdots,B$. Substituting this into \eqref{eq-S-1}, we get \eqref{eq-S-2}. This completes the proof of the corollary.
\end{proof}
\begin{rem}
Letting $k=1$ in the R.H.S. of \eqref{eq-S-1} and \eqref{eq-S-2}, we know that Corollary \ref{cor-Steklov-2} is stronger than Corollary \ref{cor-Steklov-1} by the fact that $\lambda_1(A_\Omega L_\Omega^{-1}A_B)\geq 0$.
\end{rem}
Furthermore, by applying Theorem \ref{thm-Weyl} to
$$\Delta^N u=\Delta_\Omega u+\Deg_B\cdot u-A_BL_B^{-1}A_\Omega u$$
in Proposition \ref{prop-N-Lap}, we have the following comparison on Neumann eigenvalues and Laplacian eigenvalues of the interior subgraph.
\begin{thm}\label{thm-N-IL}
Let $(G,B,m,w)$ be a  weighted finite graph with boundary such that each connected component of $G$ contains some interior vertices. Then
\begin{enumerate}
\item for $i=1,2,\cdots,|\Omega|$, $j=0,1,\cdots,|\Omega|-i$,
$$\nu_{i}\leq \mu_{i+j}(\Omega)+\lambda_{|\Omega|-j}(\Deg_B\cdot I_\Omega-A_BL_B^{-1}A_\Omega)$$
with equality for some pair $i,j$ if and only if there is a nonzero function $u\in \R^\Omega$ such that $\Delta^Nu=\nu_i u$, $\Delta_\Omega u=\mu_{i+j} u$ and $$\Deg_B\cdot u-A_BL_B^{-1}A_\Omega u=\lambda_{|\Omega|-j}(\Deg_B\cdot I_\Omega-A_BL_B^{-1}A_\Omega)u;$$
\item for $i=1,2,\cdots,|\Omega|$, $j=1,2,\cdots,i$,
$$\nu_i\geq \mu_{i-j+1}(\Omega)+\lambda_j(\Deg_B\cdot I_\Omega- A_BL_B^{-1}A_\Omega)$$
with equality for some pair $i,j$ if and only if there is a nonzero function $u\in \R^\Omega$ such that $\Delta^Nu=\lambda_i u$, $\Delta_\Omega u=\mu_{i-j+1} u$ and $$\Deg_B\cdot u-A_BL_B^{-1}A_\Omega u=\lambda_{j}(A_BL_B^{-1}A_\Omega)u.$$
\end{enumerate}
\end{thm}
\begin{rem}
By Corollary \ref{cor-Deg}, $\Deg_B\cdot I_\Omega-A_BL_B^{-1}A_\Omega$ is nonnegative. So, Theorem \ref{thm-N-IL} is stronger than Theorem \ref{thm-N-L}.
\end{rem}
Finally, translating the last result into Steklov eigenvalues, we can  obtain the following estimate on Steklov eigenvalues which is stronger than Corollary \ref{cor-Steklov-0}.
\begin{cor}\label{cor-Steklov-3}
Let $(G,B,m,w)$ be a  weighted finite graph with boundary such that each connected component of $G$ contains some boundary vertices. Then
\begin{equation*}
\begin{split}
\mu_{i-j+1}(B)+\lambda_j(\Deg_\Omega\cdot I_B- A_\Omega L_\Omega^{-1}A_B)\leq &\sigma_i(G,B)\\
\leq&\mu_{i+k}(B)+\lambda_{|B|-k}(\Deg_\Omega\cdot I_B-A_\Omega L_\Omega^{-1}A_B)
\end{split}
\end{equation*}
for $i=1,2,\cdots,|B|$, $j=1,2,\cdots,i$ and $k=0,1,2,\cdots,|B|-i$. Here $L_\Omega=\Delta^D$.
\end{cor}
\section{Some applications}
In this section, by combining the comparisons of eigenvalues in the last two sections and some known lower bounds on Laplacian eigenvalues, we get some interesting lower bounds for Neumann eigenvalues and Dirichlet eigenvalues.

First recall the Licherowicz estimates of Laplacian eigenvalues for graphs with positive curvature lower bounds.
\begin{thm}[\cite{BC,KKR, LLY}]\label{thm-Lap-Lich}
Let $(G,m,w)$ be a connected weighted finite graph with boundary.
\begin{enumerate}
\item If $(G,m,w)$ satisfies the Bakry-\'Emery curvature-dimension condition $\CD(K,n)$ with $K>0$ and $n>1$, then $\mu_2\geq \frac{nK}{n-1}$.
\item If the Ollivier curvature of $(G,m,w)$ has a positive lower bound $\kappa$, then $\mu_2\geq \kappa$.
\end{enumerate}
\end{thm}
\begin{rem}
The Ollivier curvature we used here the most general one introduced in \cite{MW} extending the definition introduced by Lin-Lu-Yau \cite{LLY} to general weighted graphs.
\end{rem}

\begin{proof}[Proof of Theorem \ref{thm-Lich-BE-1}]
By combining Theorem \ref{thm-Lap-Lich} and the comparison of Neumann eigenvalue and Laplacian eigenvalues \eqref{eq-Neumann-Laplacian}, we get Theorem \ref{thm-Lich-BE-1} directly.
\end{proof}

Similarly, by combining (2) of Theorem \ref{thm-D-L}, (2) of Theorem \ref{thm-N-IL} and Theorem \ref{thm-Lap-Lich}, we have the following Lichnerowicz-type estimates for the Dirichlet  eigenvalues and Neumann  eigenvalues directly.
\begin{thm}\label{thm-Lich-BE-2}
Let $(G,B,m,w)$ be a  connected weighted finite graph with boundary.
\begin{enumerate}
\item If $(G|_\Omega, m|_\Omega, w|_\Omega)$ is connected and satisfies the Bakry-\'Emery curvature-dimension condition $\CD(K,n)$ for some $K>0$ and $n>1$, then
\begin{equation*}
\lambda_{i}\geq \frac{nK}{n-1}+\Deg_B^{(i-1)}
\end{equation*}
and
\begin{equation*}
\nu_i\geq\frac{nK}{n-1}+\lambda_{i-1}(\Deg_B\cdot I_\Omega- A_BL_B^{-1}A_\Omega)
\end{equation*}
for $i=2,3,\cdots,|\Omega|$.
\item If $G|_\Omega$ is connected and the Ollivier curvature of $(G|_\Omega, m|_\Omega, w|_\Omega)$ has a positive lower bound $\kappa$, then,
\begin{equation*}
\lambda_i\geq \kappa+\Deg_B^{(i-1)}
\end{equation*}
and
\begin{equation*}
\nu_i\geq\kappa+\lambda_{i-1}(\Deg_B\cdot I_\Omega- A_BL_B^{-1}A_\Omega)
\end{equation*}
for $i=2,3,\cdots,|\Omega|$.
\end{enumerate}
\end{thm}
\begin{proof}
For $i=2,3,\cdots,|\Omega|$, let $j=i-1$ in (2) of Theorem \ref{thm-D-L} and (2) of Theorem \ref{thm-N-IL} and use Theorem \ref{thm-Lap-Lich} for $G|_\Omega$. We get the conclusions.
\end{proof}

Next, recall the estimate of Fiedler which relate the connectivity of graphs and their first positive Laplacian eigenvalues which are called the algebraic connectivity of graphs in \cite{FI}.
\begin{thm}[Fiedler \cite{FI}]\label{thm-Fiedler}
Let $G$ be a graph on $n$ vertices equipped with the unit weight. Then,
\begin{equation}
\mu_2\geq 2e(G)\left(1-\cos\frac{\pi}{n}\right).
\end{equation}
Here $e(G)$ is the edge connectivity of $G$. That is the least number of edges need to be deleted from $G$ to make it disconnected.
\end{thm}
By combining the comparisons of eigenvalues  with Theorem \ref{thm-Fiedler}, we have the following estimates relating the Neumann eigenvalues and Dirichlet eigenvalues to edge connectivity or edge connectivity of the interior subgraph.
\begin{thm}Let $(G,B)$ be a connected finite graph with boundary equipped with the unit weight. Then, for $i=2,3,\cdots,|\Omega|$
\begin{enumerate}
\item $\nu_2\geq 2 e(G)\left(1-\cos\frac{\pi}{|V|}\right)$;
\item $\nu_i\geq 2e(\Omega)\left(1-\cos\frac{\pi}{|\Omega|}\right)+\lambda_{i-1}\left(\Deg_B\cdot I_\Omega-A_BL_B^{-1}A_\Omega\right)$;
\item $\lambda_i\geq 2e(\Omega)\left(1-\cos\frac{\pi}{|\Omega|}\right)+\Deg_B^{(i-1)}$;
\item $\lambda_i\geq 2e(G)\left(1-\cos\frac{\pi}{|V|}\right)+\lambda_{i-1}\left(A_BL_B^{-1}A_\Omega\right)$.
\end{enumerate}
Here $e(\Omega)$ is the edge connectivity of $G|_\Omega$.
\end{thm}
\begin{proof}
By Theorem \ref{thm-Neumann-Laplacian} and Theorem \ref{thm-Fiedler}, we have
$$\nu_2\geq \mu_2\geq 2e(G)\left(1-\cos\frac{\pi}{|V|}\right).$$
This gives us (1).

In (2) of Theorem \ref{thm-N-IL}, letting $j=i-1$ and using Theorem \ref{thm-Fiedler} for $G|_\Omega$, we have
\begin{equation*}
\begin{split}
\nu_i\geq& \mu_2(\Omega)+\lambda_{i-1}\left(\Deg_B\cdot I_\Omega-A_BL_B^{-1}A_\Omega\right)\\
\geq & 2e(\Omega)\left(1-\cos\frac{\pi}{|\Omega|}\right)+\lambda_{i-1}\left(\Deg_B\cdot I_\Omega-A_BL_B^{-1}A_\Omega\right).
\end{split}
\end{equation*}
This gives us (2).

Similarly, in (2) of Theorem \ref{thm-D-L}, letting $j=i-1$ and using Theorem \ref{thm-Fiedler} for $G|_\Omega$, we get (3).

Finally, in (2) of Theorem \ref{thm-N-D-sharp-2}, letting $j=i-1$ and using (1), we get (4). This completes the proof of the theorem.
\end{proof}
Finally, recall the following geometric estimate for Laplacian eigenvalues of Friedman \cite{FR}.
\begin{thm}[Friedman \cite{FR}]\label{thm-Friedman}
Let $G$ be a connected graph on $n$ vertices equipped with the unit weight. For $2\leq i\leq n$, let $k=\left\lfloor\frac{n}{i}\right\rfloor$. Then,
\begin{enumerate}
\item when $i\not|n$,
\begin{equation*}
\mu_i\geq 2\left(1-\cos\frac{\pi}{2k+1}\right);
\end{equation*}
\item when $i|n$
\begin{equation*}
\mu_i\geq \lambda_1(P_{k+1},\{0\},m,w).
\end{equation*}
Here  $P_{k+1}$ is a path:\ $0\sim 1\sim 2\sim\cdots\sim k$, $m_j=1$ for $j=0,1,\cdots, k$, and $w_{12}=w_{23}=\cdots=w_{k-1,k}=1$ and $w_{01}=\mu_i(P_i)$ with $P_i$ equipped with the unit weight.
\end{enumerate}
\end{thm}

By combining Theorem \ref{thm-Neumann-Laplacian} and Theorem \ref{thm-Friedman}, we have the following geometric lower bounds for Neumann eigenvalues.
\begin{thm}\label{thm-Friedman-nu}
Let $(G,B)$ be a connected finite graph with boundary equipped with the unit weight. For $i=2,3,\cdots,|\Omega|$, let $k=\left\lfloor\frac{|V|}{i}\right\rfloor$.
We have
\begin{enumerate}
\item when $i \not|\ |V|$,
\begin{equation*}
\nu_i\geq 2\left(1-\cos\frac{\pi}{2k+1}\right);
\end{equation*}
\item when $i\ |\ |V|$,
\begin{equation*}
\nu_i\geq\lambda_1(P_{k+1},\{0\},m,w).
\end{equation*}
\end{enumerate}
Here $(P_{k+1},\{0\},m,w)$ is the same as in Theorem \ref{thm-Friedman}.
\end{thm}
\begin{proof}
By Theorem \ref{thm-Neumann-Laplacian}, we have
$$\nu_i\geq \mu_i.$$
Then, by Theorem \ref{thm-Friedman}, we get the conclusion.
\end{proof}
By combining (2) of Theorem \ref{thm-N-D-sharp-2} and Theorem \ref{thm-Friedman-nu}, we have the following lower bounds for Dirichlet eigenvalues.
\begin{thm}
Let $(G,B)$ be a connected finite graph with boundary equipped with the unit weight. For $i=2,3,\cdots,|\Omega|$ and $j=1,2,\cdots,i-1$, let $k=\left\lfloor\frac{|V|}{i-j+1}\right\rfloor$. We have
\begin{enumerate}
\item when $(i-j+1) \not|\ |V|$,
\begin{equation*}
\lambda_i\geq 2\left(1-\cos\frac{\pi}{2k+1}\right)+\lambda_{j}(A_BL_B^{-1}A_\Omega);
\end{equation*}
\item when $(i-j+1)\ |\ |V|$,
\begin{equation*}
\lambda_i\geq \lambda_1(P_{k+1},\{0\},m,w)+\lambda_{j}(A_BL_B^{-1}A_\Omega).
\end{equation*}
 \end{enumerate}
Here $(P_{k+1},\{0\},m,w)$ is the same as in Theorem \ref{thm-Friedman}.
\end{thm}
\begin{proof}
By (2) of Theorem \ref{thm-N-D-sharp-2}, we have
$$\lambda_i\geq \nu_{i-j+1}+\lambda_j(A_BL_B^{-1}A_\Omega).$$
Then,  by Theorem \ref{thm-Friedman-nu}, we get the conclusions.
\end{proof}
Moreover, by combining (2) of Theorem \ref{thm-D-L}, (2) of Theorem \ref{thm-N-IL} and Theorem \ref{thm-Friedman} for the interior subgraph, we have the following lower bounds for Dirichlet eigenvalues and Neumann eigenvalues.
\begin{thm}
Let $(G,B)$ be a connected finite graph with boundary equipped with the unit weight such that $G|_\Omega$ is also connected. For $i=2,3,\cdots,|\Omega|$, and $j=1,2,\cdots,i-1$,  let $k=\left\lfloor\frac{|\Omega|}{i-j+1}\right\rfloor$. We have
\begin{enumerate}
\item when $(i-j+1) \not|\ |\Omega|$,
\begin{equation*}
\lambda_i\geq 2\left(1-\cos\frac{\pi}{2k+1}\right)+\Deg_B^{(j)}
\end{equation*}
and
\begin{equation*}
\nu_i\geq 2\left(1-\cos\frac{\pi}{2k+1}\right)+\lambda_j\left(\Deg_B\cdot I_\Omega-A_BL_B^{-1}A_\Omega\right);
\end{equation*}
\item when $(i-j+1)\ |\ |\Omega|$,
\begin{equation*}
\lambda_i\geq \lambda_1(P_{k+1},\{0\},m,w)+\Deg_B^{(i)}
\end{equation*}
and
\begin{equation*}
\lambda_i\geq \lambda_1(P_{k+1},\{0\},m,w)+\lambda_{j}\left(\Deg_B\cdot I_\Omega-A_BL_B^{-1}A_\Omega\right).
\end{equation*}
 \end{enumerate}
Here $(P_{k+1},\{0\},m,w)$ is the same as in Theorem \ref{thm-Friedman}.
\end{thm}
\begin{proof}
By (2) of Theorem \ref{thm-D-L} and (2) of Theorem \ref{thm-N-IL}, we have
$$\lambda_i\geq\mu_{i-j+1}(\Omega)+\Deg_B^{(j)}$$
and
$$\mu_i\geq \mu_{i-j+1}(\Omega)+\lambda_{j}\left(\Deg_B\cdot I_\Omega-A_BL_B^{-1}A_\Omega\right).$$
Then, by applying Theorem \ref{thm-Friedman} on $G|_\Omega$, we get the conclusions.
\end{proof}

\end{document}